\documentclass[12pt]{article}
\usepackage{textcomp}
\usepackage{graphicx}           % 插入图片
\usepackage{epstopdf}
\usepackage{color}              % 支持彩色
\usepackage{xcolor,graphicx}
\usepackage{geometry}
\usepackage{float}
\usepackage{indentfirst}        % 首行缩进宏包
\usepackage{amsmath,amssymb,bm} % AMS math 宏包与数学符号加粗
\usepackage{amsthm}
\usepackage{cases}              % 数学公式
\usepackage{setspace}
\usepackage{hyperref}
\usepackage{titlesec}
\usepackage[all]{xy}            % 交换图
\usepackage{tikz}
\usepackage{multirow}% 合并 Table 的单元格时用到
\usepgflibrary{arrows}
\hypersetup{dvips,ps2pdf,CJKbookmarks,hypertex,
	colorlinks=true,
	pdfpagemode=FullScreen,
	bookmarksopen=true,
	bookmarksnumbered=true}     % 在中文环境下设置超链接文本
\usepackage{graphicx} %use graph format
\usepackage{epstopdf}
\usepackage{caption}
\usepackage{diagbox}
\usepackage{mathrsfs}

\usepackage{dsfont}
\usepackage{graphics}
\usepackage{hyperref,cleveref,cite}
\usepackage{titlesec}
\usepackage{tcolorbox}
\usepackage{fancyhdr,fancyvrb}
\usepackage{xcolor,graphicx,geometry,float,indentfirst,cases,setspace,titlesec}
\usepackage{diagbox}
\usepackage{times}

\newtheorem{con}{Conjecture}[section]
\newtheorem{defi}[con]{Definition}

\newtheorem{thm}[con]{Theorem}

\begingroup
\makeatletter \@addtoreset{equation}{section} \makeatother
\makeindex 
\def\qed{\hfill \rule{4pt}{7pt}}

 \def\remark{\vskip 0.2cm {\noindent \it Remark.}\quad}

\begin{document}
\begin{center}

 {\Large \bf A Refinement of a Theorem of Diaconis-Evans-Graham}

\end{center}

\begin{center}
 Lora R. Du$^{1}$ and {Kathy Q. Ji}$^{2}$    \vskip 2mm

$^{1,\,2}$ Center for Applied Mathematics and KL-AAGDM \\[3pt]
Tianjin University\\[3pt]
Tianjin 300072, P.R. China\\[6pt]
   \vskip 2mm

 Emails: $^1$loradu@tju.edu.cn   and   $^2$kathyji@tju.edu.cn
\end{center}
\vskip 2mm
\noindent
{\bf Abstract:} The note is dedicated to refining a theorem by Diaconis, Evans, and Graham concerning successions and fixed points of permutations. This refinement specifically addresses non-adjacent successions, predecessors, excedances, and drops of permutations.

\vskip 2mm



 \vskip 2mm

\section{Introduction}

The main objective of this paper is to give a refinement of a theorem of Diaconis-Evans-Graham \cite{Diaconis-Evans-Graham-2014} on successions and fixed points of permutations.

Let $\mathfrak{S}_n$ be the set of permutations on $[n]=\{1,2,\ldots,n\}$. For a permutation $\sigma=\sigma_1\cdots \sigma_n\in \mathfrak{S}_n$,  an index $1\leq i\leq n-1$ is called a succession if $\sigma_i+1=\sigma_{i+1}$, whereas an index $1\leq i\leq n$ is called a fixed point if $\sigma_i=i$. Let  ${\rm Suc}(\sigma)$ be the set of successions of $
\sigma$, that is
\[ {\rm Suc}(\sigma) = \{1\leq i\leq n-1~|~ \sigma_i +1=\sigma_{i+1}\}\]
and let $\overline{\rm {Fix}}(\sigma)$ denote the set of fixed points of $\sigma$ distinct from $n$. To wit,
\[  \overline{\rm {Fix}}(\sigma) = \{1\leq i\leq n-1~|  \sigma_i=i\}.\]
It should be noted that the index $n$ is excluded in the definition of $\overline{\rm {Fix}}(\sigma)$.

Given a subset $I\subseteq [n-1]$, let ${\rm Suc}_n(I)$ be the set of permutations $\sigma$ of $[n]$ such that ${\rm Suc}(\sigma)=I$ and let $\overline{\rm {Fix}}_n(I)$ be the set of permutations $\sigma\in \mathfrak{S}_n$ such that $\overline{\rm {Fix}}(\sigma)=I$.

Diaconis, Evans and Graham \cite{Diaconis-Evans-Graham-2014} discovered the following beautiful result.
\begin{thm}\label{Main} {\rm (Diaconis-Evans-Graham)}
    Let $n\geq 1$ and $I\subseteq [n-1]$. Then there is a bijection between ${\rm Suc}_n(I)$ and $\overline{\rm {Fix}}_n(I)$.
\end{thm}

It is worth mentioning that Chen \cite{Chen-1996} provided a bijective proof of the Diaconis-Evans-Graham theorem for the  case   $I=\emptyset$  via the first fundamental transformation. Brenti and Marietti \cite{Brenti-Marietti-2018} extended  this result  within the context of  colored permutations in the complex reflection groups $G(r,p,n)$ where $r,p,n$ are positive integers with $p$ dividing $n$. Recently, Chen and Fu \cite{Chen-Fu-2024} established  a left succession analogue of the Diaconis-Evans-Graham theorem, exemplifying the idea of a grammar assisted bijection.   Additionally, Ma, Qi, Yeh and Yeh \cite{Ma-Qi-Yeh-Yeh-2024} utilized the grammatical labeling technique to demonstrate  that two triple set-valued statistics of permutations are quidistributed on symmetric groups.  This implies that  the number of permutations in $\mathfrak{S}_n$ with the given  set $I$ of fixed points distinct from $1$ equals to the number of permutations in $\mathfrak{S}_n$ having $I$ as a set of $\sigma_{i+1}$ such that $\sigma_i+1=\sigma_{i+1}$ for $1\leq i\leq n-1$.

Inspired by a recent work of Chen and Fu \cite{Chen-Fu-2024}, we discover  a refinement of the Diaconis-Evans-Graham theorem  involving two  variations of successions, that is, non-adjacent successions and predecessors. Recall that Diaconis, Evans, and Graham refer to a succession of $\sigma=\sigma_1\cdots \sigma_n$ as  an unseparated pair $(k, k+1)$ of  $\sigma$ provided that $\sigma_k+1=\sigma_{k+1}$. This terminology and  the motivation for studying this concept stem from regarding a permutation as the outcome of shuffling a deck of $n$ cards. The succession has also been extensively  studied in the literature, see, e.g., \cite{Brenti-Marietti-2018, Chen-Fu-2024, Ma-Qi-Yeh-Yeh-2024, Mansour-Shattuck-2016, Reilly-Tanny-1980, Roselle-1968, Tanny-1976} , and the references cited there.

 \begin{defi}[Non-adjacent succession] Given a permutation $\sigma=\sigma_1\cdots \sigma_n\in \mathfrak{S}_n$,  an index $i$ {\rm (}$1\leq i\leq n-2${\rm )} is called a  {\it non-adjacent succession} of $\sigma$ if there  exists an integer $i+2\leq j\leq n$  such that  $\sigma_{j}=\sigma_{i}+1$. The set of non-adjacent successions of $\sigma$ is denoted by ${\rm najSuc}(\sigma)$.
 \end{defi}

\begin{defi}[Predecessor]  Given a permutation $\sigma=\sigma_1\cdots \sigma_n\in \mathfrak{S}_n$, an index $i$ {\rm (}$2\leq i\leq n${\rm )} is called a  {\it  predecessor} of $\sigma$ if there  exists an integer $1\leq j< i$   such that  $\sigma_{j}=\sigma_{i}+1$.  The set of predecessors of $\sigma$ is denoted by ${\rm Pred}(\sigma)$.
\end{defi}

For the permutation $\sigma=4\,1\,2\,6\,7\,5\,3$, we see that
\[{\rm Suc}(\sigma)=\{2,4\},\quad {\rm najSuc}(\sigma)=\{1,3\}, \quad \text{and} \quad {\rm Pred}(\sigma)=\{6,7\}.\]

To state our refinement, we also need to recall  an excedance and a drop of a permutation. For a permutation $\sigma \in \mathfrak{S}_n$, an index $1\leq i\leq n$ is called an excedance if  $\sigma_i>i$ and an index $1\leq i\leq n$ is called a drop if  $\sigma_i<i$. Define
\[  \overline{\rm Drop}(\sigma) = \{\sigma_i~| ~ 1\leq i\leq n-1,  \sigma_i<i\},\]
\[  \overline{\rm Exc}(\sigma) = \{\sigma_i ~|~ 1\leq i\leq n-1, \sigma_i>i\}.\]

It should be noted that the index $n$ is excluded in the definition of $\overline{\rm Drop}(\sigma)$ and the set $\overline{\rm Exc}(\sigma)$.
We have the following result.

  \begin{thm}%\label{Main}
    For $n\geq 1$, there is a bijection $\phi$ between $\mathfrak{S}_n$ and $\mathfrak{S}_n$ such that for $\sigma \in \mathfrak{S}_n$ and  $\tau=\phi(\sigma)$, we have
    \begin{equation}\label{rel}
    \overline{\rm {Fix}}(\sigma)={\rm Suc}(\tau), \quad \overline{\rm {Drop}}(\sigma)={\rm najSuc}(\tau) \quad  \text{and} \quad   \overline{\rm Exc}(\sigma)={\rm Pred}(\tau).
    \end{equation}
\end{thm}

\proof Given a permutation $\sigma=\sigma_1\sigma_2\cdots \sigma_n \in \mathfrak{S}_n$, we define $\tau=\phi(\sigma)$ via three steps:

{\bf Step 1.} Define $\overline{\sigma}=\overline{\sigma}_1\overline{\sigma}_2\cdots \overline{\sigma}_n$, where for $1\leq i\leq n$,
\[\overline{\sigma}_i=n+1-\sigma_{n-i+1}.\]

{\bf Step 2.} Let $\hat{\sigma}=\hat{\sigma}_1\hat{\sigma}_2\cdots \hat{\sigma}_n$, where $\hat{\sigma}_i=\overline{\sigma}_{i+1}$, for $ 1\leq i\leq n-1$ and $\hat{\sigma}_n=\overline{\sigma}_1.$   Then we write   $\hat{\sigma}$ in cycle form $(a_1,a_2,\ldots, a_r)(b_1,b_2,\ldots, b_s)\cdots (c_1,c_2,\ldots, c_t)$, where
\begin{itemize}
\item  the first cycle is the cycle including $n$, where  $n$ is placed as the  last element in this cycle;

\item other cycles are written with its smallest element first and the cycles are written in decreasing order of their smallest element.

\end{itemize}
 Define $\overline{\tau}$   to be the permutation obtained from $\hat{\sigma}$  by writing it in the above cycle form and erasing the parentheses. It can be easily verified that $\hat{\sigma}$ can be uniquely reconstructed from $\overline{\tau}$. To achieve this, we begin by inserting the first left parenthesis before $\overline{\tau}_1$ and the first right parenthesis after $n$. Then, we insert a left parenthesis before each left-to-right minimum occurring after $n$ in $\overline{\tau}$. Finally, we place a right parenthesis preceding each internal left parenthesis and at the end to obtain $\hat{\sigma}$.

{\bf Step 3.} Take the inversion of $\overline{\tau}$, denoted by   $\overline{\tau}^{-1}=\overline{\tau}_1^{-1}\cdots \overline{\tau}_n^{-1}$. Define
\[\tau=\phi(\sigma)=\tau_1\cdots \tau_n, \quad \text{where } \tau_i=n+1-\overline{\tau}^{-1}_{n-i+1} \quad  \text{ for } 1\leq i\leq n.\]
  We proceed to demonstrate that $\sigma$ and $\tau=\phi(\sigma)$ satisfy the relations \eqref{rel}.

Let
  \[ k \in {\rm \overline{Fix}}(\sigma),\quad  \sigma_r\in {\rm \overline{Drop}}(\sigma),  \quad \text{and} \quad \sigma_s\in \rm {\overline{Exc }}(\sigma).
  \]
  By definition,   we see that $\sigma_k=k$, $\sigma_r<r$ and $\sigma_s>s$. Moreover, $k,r,s\neq n$.

  Set $K=n+1-k$, $R=n+1-r$ and $S=n+1-s$. Since $k,r,s\neq n$, we see that $K,R,S\neq 1$.

  From the construction of the first step of the bijection $\phi$, we see that   \[\overline{\sigma}_K=K, \quad \overline{\sigma}_{R}=n+1-\sigma_r>R, \quad \text{and} \quad \overline{\sigma}_{S}=n+1-\sigma_s<S.\]
  Moreover, according to the construction of the second step of the bijection $\phi$, we have
   \begin{equation} \label{pf-ee}
   \hat{\sigma}_{K-1}=\overline{\sigma}_K=K, \quad \hat{\sigma}_{R-1}=\overline{\sigma}_{R}>R, \quad \text{and} \quad \hat{\sigma}_{S-1}=\overline{\sigma}_{S}<S.
   \end{equation}
   If we write the cycle decomposition of $\hat{\sigma}$ in the cycle representation described above, then there will be a cycle of the form $(\ldots,K-1, K, \ldots)$.   After the parentheses are removed to form $\overline{\tau}$, we will have $\overline{\tau}_j=K-1$ and $\overline{\tau}_{j+1}=K$ for some $1\leq j\leq n-1$.  Hence $\overline{\tau}^{-1}_{K-1}=j$, $\overline{\tau}^{-1}_{K}=j+1$, and so
   \[\tau_{k}=n+1-\overline{\tau}^{-1}_{n+1-k}=n-j  \quad \text{and} \quad \tau_{k+1}=n+1-\overline{\tau}^{-1}_{n-k}=n+1-j. \]
It follows that $k\in {\rm Suc}(\tau)$.

We proceed to show that $\sigma_r  \in {\rm najSuc}(\tau) $. Similarly,  under the assumption of the cycle form, there will be a cycle of the form  $(\ldots,R-1, \overline{\sigma}_R, \ldots)$ in the cycle representation of $\hat{\sigma}$. After the parentheses are removed to form $\overline{\tau}$, we will have    $\overline{\tau}_i=R-1$ and $\overline{\tau}_{i+1}=\overline{\sigma}_R>R$  for some $1\leq i\leq n-1$.  Hence $\overline{\tau}^{-1}_{R-1}=i$, $\overline{\tau}^{-1}_{\overline{\sigma}_R}=i+1$, and so
   \[\tau_{r+1}=n+1-\overline{\tau}^{-1}_{R-1}=n+1-i \quad \text{and} \quad \tau_{n+1-\overline{\sigma}_R}=n+1-\overline{\tau}^{-1}_{\overline{\sigma}_R}=n-i. \]
Since $n+1-\overline{\sigma}_R=\sigma_r<r$, we derive that  $\sigma_r  \in {\rm najSuc}(\tau) $.

It remains to show that $\sigma_s \in {\rm Pred}(\tau)$. By \eqref{pf-ee}, we see that $\hat{\sigma}_{S-1}\leq S-1$. If we express the cycle decomposition of $\hat{\sigma}$  using the cycle representation described above,  then there will be two situations:  a cycle of the form $(\ldots, S-1, \hat{\sigma}_{S-1},\ldots)$ or a cycle of the form $(\hat{\sigma}_{S-1},\ldots, S-1)$ occurs in the cycle decomposition of $\hat{\sigma}$. In particular, if $\hat{\sigma}_{S-1}=S-1$, then there will be a $1$-cycle $(\hat{\sigma}_{S-1})$. This case can be regarded as a special case of the situation  where $(\hat{\sigma}_{S-1},\ldots, S-1)$ occurs.

(a) If  a cycle of the form $(\ldots, S-1, \hat{\sigma}_{S-1},\ldots)$ occurs in the cycle decomposition of $\hat{\sigma}$, then $\overline{\sigma}_S=\hat{\sigma}_{S-1}\leq S-2$, and so $\sigma_s\geq s+2$.  After the parentheses are removed to obtain $\overline{\tau}$, we will have    $\overline{\tau}_t=S-1$ and $\overline{\tau}_{t+1}=\hat{\sigma}_{S-1}=\overline{\sigma}_{S}\leq S-2$  for some $1\leq t\leq n-1$.  Hence $\overline{\tau}^{-1}_{S-1}=t$, $\overline{\tau}^{-1}_{\overline{\sigma}_S}=t+1$, and so
   \[\tau_{s+1}=n+1-\overline{\tau}^{-1}_{S-1}=n+1-t \quad \text{and} \quad \tau_{n+1-\overline{\sigma}_S}=n+1-\overline{\tau}^{-1}_{\overline{\sigma}_S}=n-t. \]
Since $n+1-\overline{\sigma}_S=\sigma_s>s+2$, we derive that  $\sigma_s  \in {\rm Pred}(\tau)$.

(b) If  a cycle of the form $(\hat{\sigma}_{S-1},\ldots, S-1)$ occurs in the cycle decomposition of $\hat{\sigma}$, then  the element $n$ is not in this cycle according to the cycle form described  above, and so $(\hat{\sigma}_{S-1},\ldots, S-1)$ lies after the first cycle including $n$.  Erase the parentheses   to get $\overline{\tau}$. We will have    $\overline{\tau}_{t+1}=\hat{\sigma}_{S-1}=\overline{\sigma}_{S}$  for some $1\leq t\leq n-1$.  Since the cycles except for the first cycle are written with its smallest element first and the cycles are written in decreasing order of their smallest element, we deduce that $\overline{\tau}_{t}>\overline{\tau}_{t+1}=\overline{\sigma}_{S}$. Assume that $\overline{\tau}_{t}=T$. Hence $\overline{\tau}^{-1}_{T}=t$, $\overline{\tau}^{-1}_{\overline{\sigma}_S}=t+1$, and so
   \[\tau_{n+1-T}=n+1-\overline{\tau}^{-1}_{T}=n+1-t \quad \text{and} \quad \tau_{n+1-\overline{\sigma}_S}=n+1-\overline{\tau}^{-1}_{\overline{\sigma}_S}=n-t. \]
Since $\sigma_s =n+1-\overline{\sigma}_S> n+1-T$, we derive that  $\sigma_s  \in {\rm Pred}(\tau)$.

It is straightforward to verify that this process is reversible, and the reversed process also satisfies the relations \eqref{rel}. Thus, we complete the proof of the theorem. \qed

\remark Below is an example of the construction of $\phi(\sigma)$ from the same permutation $\sigma= 7\,2\,6\,4\,1\,3\,5 $ given by Diaconis, Evans and Graham in \cite[Remark 4.2]{Diaconis-Evans-Graham-2014}.

{\bf Step 1.} We first set $\overline{\sigma}=3\,5\,7\,4\,2\,6\,1$.

{\bf Step 2.} We then define $\hat{\sigma}= 5\,7\,4\,2\,6\,1\,3$ and we adopt the following cycle form of $\hat{\sigma}$: $(3\,4\,2\,7)(1\,5\,6).$ Thus, $\overline{\tau}=3\,4\,2\,7\,1\,5\,6.$

{\bf Step 3.}  Take the inversion of $\overline{\tau}$, denoted by $\overline{\tau}^{-1}= 5\,3\,1\,2\,6\,7\,4$. Let
  \[  \tau=\phi(\sigma)= 4\,1\,2\,6\,7\,5\,3,\]
  which differs from $ \widehat{\rho}(\sigma)= 7\,1\,2\,5\,6\,4\,3$ as obtained by Diaconis, Evans and Graham \cite{Diaconis-Evans-Graham-2014} through their bijection.

It is apparent that
\[\overline{\rm {Fix}}(\sigma)={\rm Suc}(\tau)=\{2,4\},\  \overline{\rm {Drop}}(\sigma)={\rm najSuc}(\tau)=\{1,3\}, \quad   \overline{\rm Exc}(\sigma)={\rm Pred}(\tau)=\{6,7\}.\]

 \vskip 0.2cm
\noindent{\bf Acknowledgment.} This work
was supported by the National Natural Science Foundation of China.


\begin{thebibliography}{99}

%\bibitem{Adin-Roichman-2001}

\bibitem{Brenti-Marietti-2018} F. Brenti and M. Marietti, Fixed points and adjacent ascents for classical complex reflection groups, Adv. in Appl. Math. 101 (2018) 168--183.

\bibitem{Chen-1996} W.Y.C. Chen, The skew, relative, and classical derangements, Discrete Math. 160 (1996) 235--239.

\bibitem{Chen-Fu-2024} W.Y.C. Chen and A.M. Fu, A grammar of Dumont and a theorem of Diaconis-Evans-Graham, arXiv: 2402.02743.


\bibitem{Diaconis-Evans-Graham-2014} P. Diaconis, S.N. Evans and R. Graham, Unseparated pairs and fixed points in random permutations, Adv. in Appl. Math., 61 (2014) 102--124.



\bibitem{Ma-Qi-Yeh-Yeh-2024} S.-M. Ma, H. Qi, J. Yeh and Y.-N. Yeh, On the joint distributions of succession and Eulerian statistics, arXiv:  2401.01760v2.

\bibitem{Mansour-Shattuck-2016}  T. Mansour and M. Shattuck, Counting permutations by the number of succcessions within cycles, Discrete Math. 339 (2016) 1368--1376.

\bibitem{Reilly-Tanny-1980}  J. Reilly and S. Tanny, Counting permutations by successions and other figures, Discrete Math. 32 (1980) 69--76.




\bibitem{Roselle-1968} D.P. Roselle, Permutations by number of rises and successions, Proc. Amer. Math. Soc. 19 (1968) 8--16.



 \bibitem{Tanny-1976} S. Tanny, Permutations and successions, J. Combin. Theory Ser. A 21 (1976) 196--202.

\end{thebibliography}
\end{document}